\documentclass[11pt]{article}
\usepackage[matrix,arrow,curve,cmtip]{xy}
\usepackage{amssymb}
\usepackage{latexsym}
\usepackage{theorem}
\usepackage[a4paper,width=160mm,height=230mm]{geometry}

\def\to{\mbox{$\xymatrix@1@C=5mm{\ar@{->}[r]&}$}}
\def\tto{\mbox{$\xymatrix@1@C=5mm{\ar@{=>}[r]&}$}}

\newtheorem{theorem}{Theorem}[section]
\newtheorem{lemma}[theorem]{Lemma}
\newtheorem{definition}[theorem]{Definition} 
\newtheorem{proposition}[theorem]{Proposition}
{\theorembodyfont{\upshape}\newtheorem{example}[theorem]{Example}}
\newcommand{\proof}{\noindent {\em Proof\ }: }
\def\endofproof{$\mbox{ }\hfill\Box$\par\vspace{1.8mm}\noindent}

\def\Rel{{\sf Rel}}
\def\o{^{\sf o}}
\def\bigmid{\ \Big|\ }
\def\Sh{{\sf Sh}}

\def\Matr{{\sf Matr}}
\def\Mod{{\sf Mod}}
\def\Ord{{\sf Ord}}
\def\:{\colon}
\def\impl{\Rightarrow}
\def\2{{\bf 2}}
\def\Set{{\sf Set}}
\def\op{^{\sf op}}
\def\dom{{\sf dom}}
\def\Sup{{\sf Sup}}
\def\Dist{{\sf Dist}}
\def\Map{{\sf Map}}
\def\id{{\sf id}}
\def\C{{\cal C}}
\def\E{{\cal E}}
\def\Q{{\cal Q}}
\def\R{{\cal R}}
\def\G{{\cal G}}
\def\O{{\cal O}}

\def\tensor{\cdot}
\def\<{\langle}
\def\>{\rangle}
\def\Dwn{{\sf Dwn}}
\def\can{_{\sf can}}
\def\Hilb{{\sf Hilb}}
\def\Proj{{\sf Proj}}
\def\down{\begin{picture}(6,0)%
          \put(0,0){$\downarrow$}%
          \end{picture}}

\title{Modules on involutive quantales: \\
canonical Hilbert structure,
applications to sheaf theory}
\author{Hans Heymans\footnote{Department of Mathematics and Computer Science, University of Antwerp, Middelheimlaan 1, 2020 Antwerpen, Belgium, {\tt hans.heymans@ua.ac.be}} \ and Isar Stubbe\footnote{Postdoctoral Fellow of the Research Foundation Flanders (FWO), Department of Mathematics and Computer Science, University of Antwerp, Middelheimlaan 1, 2020 Antwerpen, Belgium, {\tt isar.stubbe@ua.ac.be}}}
\date{Written: August 2008 \\ Submitted: 10 September 2008 \\ Revised: 18 May 2009}

\begin{document}

\maketitle

\begin{abstract}
We explain the precise relationship between two module-theoretic descriptions of sheaves on an involutive quantale, namely the description via so-called Hilbert structures on modules and that via so-called principally generated modules. For a principally generated module satisfying a suitable symmetry condition we observe the existence of a canonical Hilbert structure. We prove that, when working over a modular quantal frame, a module bears a Hilbert structure if and only if it is principally generated and symmetric, in which case its Hilbert structure is necessarily the canonical one. We indicate applications to sheaves on locales, on quantal frames and even on sites.
\end{abstract}

\section{Introduction}\label{A}

Jan Paseka [1999, 2002, 2003] introduced the notion of {\em Hilbert module} on an involutive quantale: it is a module equipped with an {\em inner product}. This provides for an order-theoretic notion of ``inner product space'', originally intended as a generalisation of complete lattices with a duality. Recently, Pedro Resende and Elias Rodrigues [2008] applied this definition to a locale $X$ and further defined what it means for a Hilbert $X$-module to have a {\em Hilbert basis}. These Hilbert $X$-modules with Hilbert basis describe, in a module-theoretic way, the sheaves on $X$.

At the same time, the present authors defined the notion of {\em (locally) principally generated module} on a quantaloid [Heymans and Stubbe, 2009]. Our aim too was to describe ``sheaves as modules'', albeit sheaves on quantaloids in the sense of [Stubbe, 2005b]. In this formulation the ordinary sheaves on a locale $X$ are described as locally principally generated $X$-modules whose locally principal elements satisfy an extra ``openness'' condition.

Whereas Hilbert locale modules easily generalise to modules on involutive quantales, the principally generated quantaloid modules straightforwardly specialise to involutive quantales. Thus we have two module-theoretic approaches to sheaves on involutive quantales: in this note we explain the precise relationship between them. 

This work can be summarised as follows: After some preliminary definitions we show in Section \ref{B} that any principally generated module on an involutive quantale comes with a {\em canoncial (pre-)inner product}. In Section \ref{C} we first present the notion of Hilbert basis for modules on an involutive quantale [Resende, 2008]. After introducing a suitable notion of symmetry for such modules, termed {\em principal symmetry}, we prove that a module is principally generated and principally symmetric if and only if it admits a canonical Hilbert structure (= canonical inner product plus canonical Hilbert basis). When working over a {\em modular quantal frame} it is a fact, as we prove in Section \ref{D}, that a module bears a Hilbert structure if and only if it is principally generated and principally symmetric, in which case the given inner product is necessarily the canonical one (admitting the canonical Hilbert basis). That is to say, in this case the only possible (and thus the only relevant) Hilbert structure is the canonical one. We illustrate all this module-theory with many examples. In the final Section \ref{E} we draw some conclusions from our work.

We explain all new results in this paper in a self-contained manner in the language of quantale modules, focussing on the purely order-theoretic aspects. However, in some examples, particularly those concerned with sheaf theory in one way or another, we freely use material from the references without recalling much of the details. Thus, the reader who is mainly interested in order theory can safely skip those examples; but the reader who is also interested in the applications to sheaf theory will most likely have to have a quick look at the cited papers too, insofar as the notions involved are not already familiar to her or him.

\section{Canonical inner product}\label{B}

We begin by recalling some definitions. Throughout this paper, $Q=(Q,\bigvee,\circ,1)$ stands for a {\em quantale}, i.e.\ a monoid in the monoidal category $\Sup$ of complete lattices and maps that preserve arbitrary suprema. Explicitly, a quantale $Q$ consists of a complete lattice $(Q,\bigvee)$ equipped with a binary operation $Q\times Q\to Q\:(f,g)\mapsto f\circ g$ and a constant $1\in Q$ such that 
$$f\circ(g\circ h)=(f\circ g)\circ h\mbox{, \quad}1\circ f=f=f\circ 1\mbox{\quad and \quad}(\bigvee_{i\in I}f_i)\circ(\bigvee_{j\in J}g_j)=\bigvee_{i\in I}\bigvee_{j\in J}(f_i\circ g_j)$$ 
for all $f,g,h,f_i,g_j\in Q$. (Some call this a {\em unital quantale}, but since we shall not encounter ``non-unital quantales'' in this work we drop that adjective.)
\begin{definition}\label{1}
A map $Q\to Q\:f\mapsto f\o$ is an {\em involution}, and the pair $(Q,(-)\o)$ forms an {\em involutive quantale}, if it is order-preserving ($f\leq g\impl f\o\leq g\o$), involutive ($f^{\sf oo}=f$) and multiplication-reversing ($(f\circ g)\o=g\o\circ f\o$). \end{definition}
It follows that an involution is an isomorphism of complete lattices, and also unit-preserving ($1\o=1$). Most often we shall simply speak of ``an involutive quantale $Q$'' and leave it understood that the involution is written as $f\mapsto f\o$. 
\begin{definition}\label{1.0}
An element $f\in Q$ of an involutive quantale is {\em symmetric} if $f\o=f$.
\end{definition}
A symmetric idempotent element of $Q$ ($f\o=f=f\circ f$) is sometimes called a {\em projection}. 
\begin{example}\label{1.1} Among the many examples of involutive quantales, we point out some of particular interest.
\begin{enumerate}
\item\label{1.1.1} A quantale $Q$ is commutative if and only if the identity map $1_Q\:Q\to Q\:q\mapsto q$ is an involution. In particular is every {\em locale} (also called {\em frame}) $X=(X,\bigvee,\wedge,\top)$ an involutive quantale for this trivial involution.
\item\label{1.1.2} Let $S$ be a complete lattice with a {\em duality}, i.e.\ a supremum-preserving map $d\:S\to S\op$ such that $d(x)=d^*(x)$ and $d(d(x))=x$ for all $x\in S$, where $d^*$ is the right adjoint to $d$ (abbreviated as $d\dashv d^*$) in the category $\Ord$ of ordered sets and order-preserving maps, explicitly, $d^*\:S\op\to S\:y\mapsto\bigvee\{x\in S\mid d(x)\leq\op y\}$. The quantale $Q(S):=(\Sup(S,S),\bigvee,\circ,1_S)$ has a natural involution [Mulvey and Pelletier, 1992]: for $f\in Q(S)$ put $f\o:=d\op\circ (f^*)\op\circ d$ (where $f\dashv f^*$ in $\Ord$). When putting $f_{\sf o}:=d\op\circ f\op\circ d$, we have $f\o\dashv f_{\sf o}$ in $\Ord$.
\item\label{1.1.3} A {\em modular quantale} $Q$ is an involutive one which satisfies Freyd's modular law [Freyd and Scedrov, 1990]: $(p\circ q\wedge r)\leq p\circ(q\wedge p\o\circ r)$ for all $p,q,r\in Q$. We follow [Resende, 2007] in speaking of a {\em quantal frame} when we mean a quantale whose underlying lattice is a frame (= locale); the term {\em modular quantal frame} then speaks for itself. It is a matter of fact that modular quantal frames are precisely the one-object locally complete distributive allegories of P. Freyd and A. Scedrov [1990]. Allegories are closely related to toposes; below we shall see that modular quantal frames in particular appear in the study of sheaves (cf.\ Theorem \ref{14} and Example \ref{20}).
\item\label{1.1.4} An {\em inverse quantal frame} is a modular quantal frame $Q$ in which every element is the join of so-called {\em partial units} ($p\in Q$ is a partial unit if $p\o p\vee pp\o\leq 1_Q$); it suffices that the top of $Q$ is such a join. This definition is equivalent to the original one given in [Resende, 2007] because it is proved in that reference that inverse quantal frames arise as quotients (as frames {\em and} as involutive quantales) of quantal frames that are evidently modular. There is a correspondence up to isomorphism between inverse quantal frames and \'etale groupoids [Resende, 2007], providing a context to consider \'etendues in terms of quantales.
\item\label{1.1.5} (In this and the next example we use notions that, for a lack of space, we cannot recall; but we do include ample references.) A {\em quantaloid} is a $\Sup$-enriched category. If $A$ is an object of a quantaloid $\Q$, then $\Q(A,A)$ is a quantale; in particular, a quantaloid with only one object is precisely a quantale. A quantaloid $\Q$ has a direct-sum completion, which can be described as $\Matr(\Q)$, the quantaloid of {\em matrices with elements in $\Q$}. All definitions above can straightforwardly be generalised from quantales to quantaloids. For details, see e.g.\ [Freyd and Scedrov, 1990; Rosenthal, 1996; Stubbe, 2005a]. A small quantaloid $\Q$ is {\em Morita equivalent} to the quantale $Q:=\Matr(\Q)(\Q_0,\Q_0)$ [Mesablishvili, 2004], and it is easily seen that several properties of $\Q$ are carried over to its Morita-equivalent quantale $Q$: for example, if $\Q$ is involutive then so is $Q$. 
\item\label{1.1.6} For a small site $(\C,J)$, i.e.\ $\C$ a small category and $J$ a Grothendieck topology on $\C$, the $J$-closed relations between the representables in $\Set^{\C\op}$ form a locally complete distributive allegory, i.e.\ a modular quantaloid $\Q$ whose hom-objects are frames [Walters, 1982; Betti and Carboni, 1983]. It is easy to verify that this small quantaloid's Morita-equivalent quantale $Q$ is a modular quantal frame, and that $\Q$ can be identified with a subquantaloid of the universal splitting of the symmetric idempotents of $Q$. 
\end{enumerate}
\end{example}

When we speak of a {\em (right) $Q$-module $M$} we mean so in the obvious way in $\Sup$. That is to say, $(M,\bigvee)$ is a complete lattice on which $Q$ acts by means of a function $M\times Q\to M\:(m,f)\mapsto m\tensor f$ satisfying
$$m\tensor(f\circ g)=(m\tensor f)\tensor g\mbox{, \quad}m\tensor 1=m\mbox{\quad and \quad}(\bigvee_{i\in I}m_i)\tensor(\bigvee_{j\in j}f_j)=\bigvee_{i\in I}\bigvee_{j\in J}(m_i\tensor f_j)$$
for all $m,m_i\in M$ and $f,g,f_j\in Q$. Accordingly, a function $\phi\:M\to N$ between two $Q$-modules is a {\em $Q$-module morphism} if 
$$\phi(m\tensor f)=\phi(m)\tensor f\mbox{\quad and \quad}\phi(\bigvee_{i\in I}m_i)=\bigvee_{i\in I}\phi(m_i)$$
for all $m,m_i\in M$ and $f\in Q$. We shall write $\Mod(Q)$ for the category of $Q$-modules and module morphisms. Of course $Q$ itself is a $Q$-module, with action given by multiplication in $Q$.
\begin{definition}[Paseka, 1999]\label{2}
Let $M$ be a module on an involutive quantale $Q$. A map
$$M\times M\to Q\:(m,n)\mapsto\<m,n\>$$
is a {\em pre-inner product} if, for all $m,n\in M$,
\begin{enumerate}
\item\label{2.1} $\<m,-\>\:M\to Q$ is a module morphism,
\item\label{2.2} $\<m,n\>\o=\<n,m\>$ (which we refer to as {\em Hermitian symmetry}).
\newcounter{saveenum}
\setcounter{saveenum}{\value{enumi}}
\end{enumerate}
It is an {\em inner product} if it moreover satisfies
\begin{enumerate}
\setcounter{enumi}{\value{saveenum}}
\item\label{2.3} $\<-,m\>=\<-,n\>$ implies $m=n$
\setcounter{saveenum}{\value{enumi}}
\end{enumerate}
and it is said to be {\em strict} if
\begin{enumerate}
\setcounter{enumi}{\value{saveenum}}
\item\label{2.4} $\<m,m\>=0$ implies $m=0$.
\end{enumerate}
\end{definition}

Now we shall recall some definitions from [Heymans and Stubbe, 2009], where they are given for quantaloids but which we apply here to quantales.

Let $e\in Q$ be an idempotent. The fixpoints of $e\circ -\:Q\to Q$ form a $Q$-module which we shall write as $Q^e$: the action of $Q$ on $Q^e$ is given by multiplication, so the inclusion
$$\iota_e\:Q^e\to Q\:f\mapsto f$$
is a module morphism. Further, if $M$ is any $Q$-module then for any $m\in M$ the map
$$\tau_m\:Q\to M\:f\mapsto m\tensor f$$
is a module morphism. Thus also the composite
$$\xymatrix{
Q^e\ar[rd]_{\iota_e}\ar@{.>}[rr]^{\zeta_m} & & M \\
 & Q\ar[ru]_{\tau_m}}$$
is a module morphism. Essentially as an application of the Yoneda Lemma for enriched categories [Kelly, 1982] we find the following characterisation.
\begin{proposition}\label{3.0}
Let $Q$ be a quantale, $e\in Q$ an idempotent, and $M$ a $Q$-module. There is a one-one correspondence between the fixpoints of $-\tensor e\:M\to M$ and the module morphisms from $Q^e$ to $M$.
\end{proposition}
\proof
If $\zeta\:Q^e\to M$ is any module morphism, then $m_{\zeta}:=\zeta(e)\in M$ satisfies $m_{\zeta}\tensor e=m_{\zeta}$; conversely, if $m\in M$ satisfies $m\tensor e=m$, then $\zeta_m\:Q^e\to M\:f\mapsto m\tensor f$ is a module morphism. This is easily seen to set up a one-one correspondence. 
\endofproof
In particular, such a map $\zeta_m\:Q^e\to M$ between complete lattices preserves suprema; therefore it has an infima-preserving right adjoint in the category of ordered sets and order-preserving maps.
However, in general the order-preserving right adjoint need not be a module morphism, i.e.\ it need not be right adjoint to $\zeta_m$ in the category $\Mod(Q)$ of $Q$-modules.
\begin{definition}[Heymans and Stubbe, 2009]\label{3}
Let $Q$ be a quantale and $M$ a $Q$-module. An element $m\in M$ is said to be {\em locally principal at an idempotent $e\in Q$} if $m\tensor e=m$ and
$$\zeta_m\:Q^e\to M\:f\mapsto m\tensor f$$
has a right adjoint in $\Mod(Q)$.
\end{definition}
\begin{proposition}\label{4}
Let $Q$ be a quantale, $e\in Q$ an idempotent, and $M$ a $Q$-module. There is a one-one correspondence between $M$'s locally principal elements at $e$ and left adjoint module morphisms from $Q^e$ to $M$.
\end{proposition}
In what follows we shall always write $\zeta^*\:M\to Q^e$ for the right adjoint to a given module morphism $\zeta\:Q^e\to M$, whenever we know or assume it exists.

Now we come to a trivial but crucial observation.
\begin{proposition}\label{6.1}
Let $Q$ be an involutive quantale, $\E\subseteq Q$ the set of symmetric idempotents, and $M$ a $Q$-module. The formula
$$\<m,n\>\can:=\bigvee\left\{(\zeta^*(m))\o\circ\zeta^*(n)\bigmid e\in\E,\ \zeta\:Q^e\to M\mbox{ left adjoint in }\Mod(Q)\right\}$$
defines a pre-inner product, called the {\em canonical pre-inner product}, on $M$.
\end{proposition}
\proof
For any $e\in\E$ and any left adjoint $\zeta\:Q^e\to M$, the pointwise multiplication of the composite module morphism
$$\xymatrix{M\ar[r]^{\zeta^*} & Q^e\ar[r]^{\iota_e} & Q}$$
with the element $(\zeta^*(m))\o\in Q$ gives a module morphism
$$(\zeta^*(m))\o\circ\zeta^*(-)\:M\to Q.$$
As any pointwise supremum of parallel module morphisms is again a module morphism, we find that 
$$\<m,-\>\can=\bigvee\left\{(\zeta^*(m))\o\circ\zeta^*(-)\bigmid e\in\E,\ \zeta\:Q^e\to M\mbox{ left adjoint in }\Mod(Q)\right\}$$
is a module morphism from $M$ to $Q$. It is a triviality that the function $\<-,-\>\can\:M\times M\to Q$ is symmetric.
\endofproof

\begin{example}\label{6.2} We shall compute some more explicit examples in the next section, but we already include the following here.
\begin{enumerate}
\item\label{6.2.1} Every involutive quantale $Q$, regarded as a module over itself, has a natural inner product [Paseka, 1999]: for $f,g\in Q$ let $\<f,g\>:=f\o\circ g$. And the canonical pre-inner product on a $Q$-module $M$ is expressed as a supremum of values of the natural inner product on $Q$.
\item\label{6.2.2} Particularly for a complete lattice $S$ with duality $d\:S\to S\op$ we can consider the natural inner product on the involutive quantale $Q(S)$: we have for $s,t\in S$ and $f,g\in Q(S)$ that
$$\<f,g\>(s)\leq t\iff g(s)\leq f_{\sf o}(t).$$ 
From this it is easy to verify that $\<f,g\>=0$ if and only if $f$ and $g$ are {\em disjoint}: $f(s)\leq d(g(t))$ for all $s,t\in S$.
\end{enumerate}
\end{example}

\section{Canonical Hilbert basis}\label{C}

We start by recalling another definition from [Heymans and Stubbe, 2009] (where it was actually stated more generally for modules on quantaloids).
\begin{definition}\label{6}
Let $Q$ be a quantale, $\E\subseteq Q$ any set of idempotent elements containing the unit $1$, and $M$ a $Q$-module. If, for all $m\in M$, 
$$m=\bigvee\left\{\zeta(\zeta^*(m))\bigmid e\in\E,\ \zeta\:Q^e\to M\mbox{ left adjoint in }
\Mod(Q)\right\}$$
then $M$ is {\em $\E$-principally generated} (which is short for: {\em generated by its elements which are locally principal at some $e\in\E$}).
\end{definition}
This Definition \ref{6} resembles the following notion, which was originally defined by Resende and Rodrigues [2008] for pre-Hilbert modules on a locale, but which can straightforwardly be extrapolated to pre-Hilbert modules on an involutive quantale, as [Resende, 2008] does:
\begin{definition}\label{7}
Let $Q$ be an involutive quantale, and $M$ a $Q$-module with pre-inner product $\<-,-\>$. If a subset $\Gamma\subseteq M$ satisfies, for all $m\in M$,
$$m=\bigvee_{s\in\Gamma}s\tensor\<s,m\>$$
then it is a {\em Hilbert basis}\footnote{As also remarked in [Resende and Rodrigues, 2008], the word ``basis'' is quite deceiving: since there is no freeness condition, it would be more appropriate to speak of {\em Hilbert generators}. However, for the sake of clarity we shall adopt the terminology that was introduced in the cited references.} for $M$.
\end{definition}
If a $Q$-module $M$ bears a pre-inner product admitting a Hilbert basis, we speak of its {\em Hilbert structure}; unless explicitly stated otherwise we shall always write $\<-,-\>$ for the pre-inner product and $\Gamma$ for the Hilbert basis. As also pointed out in [Resende, 2008], it is trivial to check that:
\begin{proposition}\label{7.1}
If $Q$ is an involutive quantale, and $M$ a $Q$-module with a pre-inner product $\<-,-\>$ admitting a Hilbert basis $\Gamma$, then $\<-,-\>$ is in fact an inner product. 
\end{proposition}
\proof
Given $m,n\in M$ such that $\<-,m\>=\<-,n\>$, we certainly have $\<s,m\>=\<s,n\>$ for all $s\in\Gamma$. The formula in Definition \ref{7} then allows us to compute that 
$$m=\bigvee_{s\in\Gamma}s\tensor\<s,m\>=\bigvee_{s\in\Gamma}s\tensor\<s,n\>=n$$
and we are done.
\endofproof
Both Definitions \ref{6} and \ref{7} speak of a ``generating set'' for $Q$-modules... But already in the localic case these two definitions are different!
\begin{example}\label{9}
Let $X$ be a locale, view it as a quantale $(X,\bigvee,\wedge,\top)$ with identity involution. The set $\E$ of symmetric idempotents in $X$ coincides with $X$, and it is shown in [Stubbe, 2005b; Heymans and Stubbe, 2009] that an $\E$-principally generated $X$-module is the same thing as an {\em ordered sheaf} on $X$, i.e.\ an ordered object in $\Sh(X)$. On the other hand, as proved in [Resende and Rodrigues, 2008], a pre-Hilbert $X$-module with Hilbert basis is the same thing as a {\em sheaf} on $X$.
\end{example}
This example hints at the importance of the intrinsic {\em symmetry} in the notion of ``pre-Hilbert $Q$-module with Hilbert basis'', i.e.\ the symmetry of the involved pre-inner product. Indeed notice that Definitions \ref{2} and \ref{7} ask for a module on an {\em involutive} quantale -- without which it would simply be impossible to coherently speak of {\em symmetry} -- whereas Definition \ref{6} has no such requirement at all. To systematically explain the relation between the two definitions we must therefore develop a suitable notion of symmetry in the context of $\E$-principally generated $Q$-modules.
\begin{proposition}\label{10}
Let $Q$ be an involutive quantale, $\E\subseteq Q$ the set of symmetric idempotents, and $M$ a $Q$-module. The following statements are equivalent:
\begin{enumerate}
\item for any $e\in\E$, any left adjoint $\zeta\:Q^e\to M$ and any $m\in M$: $\zeta^*(m)=\<\zeta(e),m\>\can$,
\item for any $e,f\in\E$ and any left adjoints $\zeta\:Q^e\to M$, $\eta\:Q^f\to M$: $\zeta^*(\eta(f))=\<\zeta(e),\eta(f)\>\can$,
\item for any $e,f\in\E$ and any left adjoints $\zeta\:Q^e\to M$, $\eta\:Q^f\to M$: $\zeta^*(\eta(f))=(\eta^*(\zeta(e)))\o$.
\end{enumerate}
In this case we say that $M$ is {\em $\E$-principally symmetric}.
\end{proposition}
\proof
The only non-trivial implication is ($3\impl 1$). In fact, the ``$\leq$'' in statement 1 always holds: because 
$$\zeta^*(m)=e\circ\zeta^*(m)=e\o\circ\zeta^*(m)\leq(\zeta^*(\zeta(e)))\o\circ\zeta^*(m)\leq\<\zeta(e),m\>\can$$
where we used respectively: $\zeta^*(m)\in Q^e$; $e=e\o$; the unit of the adjunction $\zeta\dashv\zeta^*$ to get $e\leq\zeta^*(\zeta(e))$ from which $e\o\leq(\zeta^*(\zeta(e)))\o$ because the involution preserves order; and finally the definition of the canonical pre-inner product.

Thus, assuming statement 3 we must show that the ``$\geq$'' in statement 1 holds. But we can compute that, for any $f\in\E$ and any left adjoint $\eta\:Q^f\to M$, 
$$(\eta^*(\zeta(e)))\o\circ\eta^*(m)=\zeta^*(\eta(f))\circ\eta^*(m)=\zeta^*\circ\eta\circ\eta^*(m)\leq\zeta^*(m)$$
using respectively: the assumption; the fact that $\zeta^*(\eta(f))$ is the representing element for the $Q$-module morphism $\zeta^*\circ\eta\:Q^f\to Q^e$ (cf.\ Proposition \ref{3.0}); and the counit of the adjunction $\eta\dashv\eta^*$.
\endofproof
Remark that $(1\impl 2\impl 3)$ in Proposition \ref{10} holds {\em for any pre-inner product} on $M$ (but $(3\impl 1)$ does not!): that is to say, if one can prove the first or the second condition for a given pre-inner product on $M$ (not necessarily the canonical one), then it follows that $M$ is $\E$-principally symmetric. This shall be useful in the proof of Lemma \ref{18}.

We can now prove a first ``comparison'' between Definitions \ref{6} and \ref{7}.
\begin{theorem}\label{11}
Let $Q$ be an involutive quantale, $\E\subseteq Q$ the set of symmetric idempotents, and $M$ a $Q$-module. The following are equivalent:
\begin{enumerate}
\item $M$ is $\E$-principally generated and $\E$-principally symmetric,
\item the set
$$\Gamma\can:=\{\mbox{all elements of $M$ which are locally principal at some }e\in\E\}$$
is a Hilbert basis for the canonical pre-inner product on $M$, called the {\em canonical Hilbert basis}.
\end{enumerate}
In this case, it follows by Proposition \ref{7.1} that the canonical pre-inner product is an inner product; we speak of the {\em canonical Hilbert structure on $M$}.
\end{theorem}
\proof
(1$\impl$2) Assuming that $M$ is $\E$-principally generated we have by definition that, for any $m\in M$,
$$m=\bigvee\left\{\zeta(\zeta^*(m))\bigmid e\in\E,\ \zeta\:Q^e\to M\mbox{ left adjoint in }
\Mod(Q)\right\}.$$
Assuming moreover that $M$ is $\E$-principally symmetric we can compute
$$\zeta(\zeta^*(m))=\zeta(e\circ\zeta^*(m))=\zeta(e\circ\<\zeta(e),m\>\can)=\zeta(e)\tensor\<\zeta(e),m\>\can$$
using respectively: $\zeta^*(m)\in Q^e$; the first statement in Proposition \ref{10}; and the fact that $\zeta$ is a module morphism. Replacing this in the right hand side of the first expression, we obtain
$$m=\bigvee\left\{\zeta(e)\tensor\<\zeta(e),m\>\can\bigmid e\in\E,\ \zeta\:Q^e\to M\mbox{ left adjoint in }
\Mod(Q)\right\}$$
so that, if we put
$$\Gamma\can:=\left\{\zeta(e)\bigmid e\in\E,\ \zeta\:Q^e\to M\mbox{ left adjoint in }
\Mod(Q)\right\},$$
which we know by Proposition \ref{4} indeed corresponds to the set of elements of $M$ which are locally principal at some $e\in\E$, we find precisely what we claimed.

(2$\impl$1) For any $e\in\E$ and left adjoint $\zeta\:Q^e\to M$, there certainly is a module morphism
$\<\zeta(e),-\>\can\:M\to Q$. But we can compute that, for any $m\in M$, 
$$e\circ\<\zeta(e),m\>\can=\<\zeta(e)\tensor e\o,m\>\can=\<\zeta(e\circ e\o),m\>\can=\<\zeta(e),m\>\can$$
using: the ``conjugate-linearity'' of $\<-,m\>\can$; the module morphism $\zeta$; the fact that $e$ is a symmetric idempotent. Therefore, this module morphism corestricts to $\<\zeta(e),-\>\can\:M\to Q^e$. We claim that it is right adjoint to $\zeta\:Q\to M$. Indeed, if for $q\in Q^e$ and $m\in M$ we assume that $q\leq\<\zeta(e),m\>\can$ then we can compute
$$\zeta(q)=\zeta(e\circ q)=\zeta(e)\tensor q\leq\zeta(e)\tensor\<\zeta(e),m\>\can\leq m$$
using: $q\in Q^e$, i.e.\ $e\circ q=q$; $\zeta$ is a module morphism; the assumed inequality which is preserved by $\zeta(e)\tensor -$; and finally the hypothesis that $M$ has Hilbert basis $\Gamma$. Assuming conversely that $\zeta(q)\leq m$, then we can compute
$$\begin{array}{rl}
q & =e\circ q\leq\<\zeta(e),\zeta(e)\>\can\circ q=\<\zeta(e),\zeta(e)\tensor q\>\can \\[2ex]
  & =\<\zeta(e),\zeta(e\circ q)\>\can=\<\zeta(e),\zeta(q)\>\can\leq\<\zeta(e),m\>\can
\end{array}$$
using: $q\in Q^e$; the unit of $\zeta\dashv\zeta^*$ in $e=e\o\circ e\leq(\zeta^*(\zeta(e)))\o\circ\zeta^*(\zeta(e))\leq\<\zeta(e),\zeta(e)\>\can$; the module morphism $\<\zeta(e),-\>\can$; the module morphism $\zeta$; again $q\in Q^e$; and finally the assumed inequality which is preserved by the module morphism $\<\zeta(e),-\>\can$. Hence, for any $q\in Q^e$ and $m\in M$,
$$q\leq\<\zeta(e),m\>\can\iff \zeta(q)\leq m.$$
Adjoints are unique and so we obtain that $\zeta^*(m)=\<\zeta(e),m\>\can$ for all $m\in M$. By Proposition \ref{10} this exactly means that $M$ is $\E$-principally symmetric. Since we assume that $\Gamma$ is a Hilbert basis, we have that
$$m=\bigvee\left\{\zeta(e)\tensor\<\zeta(e),m\>\can\bigmid e\in\E,\zeta\:Q^e\to M\mbox{ left adjoint in }\Mod(Q)\right\}.$$
But the previous computation allows us to write
$$\zeta(e)\tensor\<\zeta(e),m\>\can=\zeta(e)\tensor\zeta^*(m)=\zeta(e\circ\zeta^*(m))=\zeta(\zeta^*(m))$$
hence we find that
$$m=\bigvee\left\{\zeta(\zeta^*(m))\bigmid e\in\E,\zeta\:Q^e\to M\mbox{ left adjoint in }\Mod(Q)\right\}$$
as wanted.
\endofproof

\begin{example}\label{12} We shall give some examples of $Q$-modules with Hilbert structure, and then make a comment on the {\em category} of $Q$-modules with Hilbert structure.
\begin{enumerate}
\item\label{12.1} Cf.\ Example \ref{6.2}--\ref{6.2.1}, $\Gamma:=\{1_Q\}$ is a Hilbert basis for the natural inner product on $Q$. More generally, if $e\in Q$ is an idempotent, then $Q^e$ is a $Q$-module with inner product $\<f,g\>:=f\o\circ g$ admitting $\Gamma:=\{e\}$ as Hilbert basis.
\item\label{12.2} Let $Q$ be the 2-element chain $\2=\{0<1\}$ (with $\wedge$ as multiplication, trivial involution, etc.); both its elements are symmetric idempotents. Let $(A,\leq)$ be an ordered set and consider $\Dwn(A,\subseteq)$, the downclosed subsets of $A$ ordered by inclusion. This is the typical example of an $\E$-principally generated $\2$-module [Heymans and Stubbe, 2009] and is also one of the fundamental constructions in [Resende and Rodrigues, 2008]. If $D\in\Dwn(A, \subseteq)$ is a locally principal element, then it is either the empty downset $D=\emptyset$ (locally principal at $0\in\2$) or a principal downset $D=\down x$ for some $x\in A$ (locally principal at $1\in\2$). For any $D,E\in\Dwn(A,\subseteq)$, their canonical inner product is 
$$\<D,E\>\can=\left\{\begin{array}{l}1\mbox{ if }D\cap E\neq \emptyset \\ 0\mbox{ otherwise}\end{array}\right.$$
To say that $\Dwn(A,\subseteq)$ is $\E$-principally symmetric is to require that for any $x,y\in A$:
$$\down x\subseteq \down y\iff \down y\subseteq\down x.$$
This makes the order $(A,\subseteq)$ in reality an equivalence relation $(A,\approx)$.
\item\label{12.3} The localic case: Let $X$ be any locale and $S$ any set. Then $X^S$ is an $X$-module, with pointwise suprema and $(f\tensor x)(s)=f(s)\wedge x$, for any $f\in X^S$, $x\in X$ and $s\in S$. Take now an $X$-matrix $\Sigma\:S\to S$ (= a family $(\Sigma(y,x))_{(x,y)\in S\times S}$ of elements of $X$) satisfying
$$\Sigma(z,y)\wedge\Sigma(y,x)\leq\Sigma(z,x)\mbox{ \ \ and \ \ }\Sigma(x,x)\wedge\Sigma(x,y)=\Sigma(x,y)=\Sigma(x,y)\wedge\Sigma(y,y)$$
and consider the $X$-submodule $\R(\Sigma)$ of $X^S$ consisting of those functions $f\:S\to X$ satisfying 
$$f(s)=\bigvee_{x\in S}\Sigma(s,x)\wedge f(x).$$ 
In the terminology of [Stubbe, 2005b], $\Sigma$ is a totally regular $X$-semicategory and $\R(\Sigma)$ is (up to the identification of $X$-modules with cocomplete $X$-categories [Stubbe, 2006]) the cocomplete $X$-category of (totally) regular presheaves on $\Sigma$. This is the typical example of a {\em locally principally generated $X$-module} [Heymans and Stubbe, 2009] and is one of the fundamental constructions of [Resende and Rodrigues, 2008] too. It is not too difficult to show by direct calculations, but it also follows from our further results, that $\R(\Sigma)$ is $\E$-principally symmetric if and only if $\Sigma$ is a symmetric $X$-matrix. Moreover, for a symmetric $X$-matrix $\Sigma$ to satisfy the above conditions is equivalent to it being an idempotent, hence the module $\R(\Sigma)$ is $\E$-principally generated and $\E$-principally symmetric if and only if $\Sigma$ is a so-called {\em projection matrix} (with elements in $X$). Our upcoming Theorem \ref{14} says that such structures coincide in turn with $X$-modules with (necessarily canonical) Hilbert structure.
\item\label{12.4} The previous example is an instance of a more general situation. We write $\Hilb(Q)$ for the quantaloid whose objects are $Q$-modules with Hilbert structure and whose morphisms are module morphisms. And we write $\Matr(Q)$ for the quantaloid whose objects are sets and whose morphisms are {\em matrices with elements in $A$}: such a matrix $\Lambda\:S\to T$ is an indexed set of elements of $Q$, $(\Lambda(y,x))_{(x,y)\in S\times T}\in Q$. Matrices compose straightforwardly with a ``linear algebra formula'', and the identity matrix on a set $S$ has all $1$'s on the diagonal and $0$'s elsewhere. This matrix construction makes sense for any quantale (and even quantaloid), and whenever $Q$ is involutive then so is $\Matr(Q)$: the involute of a matrix is computed elementwise. Now there is an equivalence of quantaloids\footnote{[Resende, 2008] also notes the object correspondence, but not the morphism correspondence, and thus not the equivalence of these quantaloids.} 
$$\Hilb(Q)\simeq\Proj(Q),$$ 
where the latter is the quantaloid obtained by splitting the symmetric idempotents in $\Matr(Q)$, i.e.\ the quantaloid of so-called {\em projection matrices} with elements in $Q$. Explicitly, if $\Sigma\:S\to S$ is such a projection matrix, then 
$$\R(\Sigma):=\{f\:S\to Q\mid \forall s\in S:f(s)=\bigvee_{s\in S}\Sigma(s,x)\circ f(x)\}$$
is a $Q$-module with inner product and Hilbert basis respectively
$$\<f,g\>:=\bigvee_{s\in S}(f(s))\o\circ g(s)\mbox{ and }\Gamma:=\{f_s\:S\to Q\:x\mapsto\Sigma(x,s)\mid s\in S\}.$$
This object correspondence $\Sigma\mapsto\R(\Sigma)$ extends to a $\Sup$-functor from $\Proj(Q)$ to $\Hilb(Q)$: it is the restriction to symmetric idempotent matrices of the embedding of the Cauchy completion of $Q$ qua one-object $\Sup$-category -- i.e.\ the quantaloid obtained by splitting {\em all} idempotents of $\Matr(Q)$ -- into $\Mod(Q)$. Conversely, a module $M$ with inner product $\<-,-\>$ and Hilbert basis $\Gamma$ obviously determines a projection matrix $\Sigma\:\Gamma\to\Gamma$ with elements $\Sigma(s,t):=\<s,t\>$; this easily extends to a $\Sup$-functor from $\Hilb(Q)$ to $\Proj(Q)$. These two functors set up the equivalence.
\item\label{12.5} A notable consequence of the previous example is the existence of an involution on $\Hilb(Q)$, induced by the obvious involution on $\Proj(Q)$: the involute of a morphism $\phi\:M\to N$ in $\Hilb(Q)$ is the unique module morphism $\phi\o\:N\to M$ characterised by 
$$\<\phi(s),t\>=\<s,\phi\o(t)\>$$
for all basis elements $s$ of $M$ and $t$ of $N$.
\end{enumerate}
\end{example}

In the localic case (cf.\ the example above) we can moreover prove an alternative formulation of the symmetry condition in Proposition \ref{10}: an ``openness'' condition formulated on the principal elements. In the next example we recall and explain this.
\begin{example}\label{13}
Let $X$ be a locale. Every element $u\in X$ is a symmetric idempotent, and the open sublocale $\down u\subseteq X$ is precisely the $X$-module of fixpoints of $u\wedge-\:X\to X$. If $M$ is an $X$-module for which each left adjoint $\zeta\:\down u\to M$ is {\em open}, in the sense that
$$\mbox{for all $x\leq u$ and $m\in M$: }\zeta(x\wedge\zeta^*(m))=\zeta(x)\wedge m,$$
then it is $\E$-principally symmetric. The converse also holds, provided that $M$ is $\E$-principally generated, in which case $M$ is an {\em \'etale $X$-module} in the terminology of [Heymans and Stubbe, 2009].
\end{example}
\proof 
Let $\zeta\:\down u\to M$ and $\eta\:\down v\to M$ be left adjoints, suppose that $\zeta$ is open: with $x:=u$ and $m:=\eta(v)$ in the above formula, it follows that
$$\zeta(\zeta^*(\eta(v)))=\zeta(u)\wedge\eta(v).$$
Applying $\eta^*$ (which preserves infima) gives
$$(\eta^*\circ\zeta)(\zeta^*(\eta(v)))=\eta^*(\zeta(u))\wedge\eta^*(\eta(v)).$$
The right hand side equals $\eta^*(\zeta(u))$ because $\eta^*(\eta(v))=v$ (the adjunction $\eta\dashv\eta^*$ splits). The left hand side equals $\eta^*(\zeta(u))\wedge\zeta^*(\eta(v))$, because the $X$-module morphism $\eta^*\circ\zeta\:\down u\to\down v$ is represented by $\eta^*(\zeta(u))\leq u\wedge v$. Thus we get 
$$\eta^*(\zeta(u))\wedge\zeta^*(\eta(v))=\eta^*(\zeta(u))\mbox{, or in other words }\eta^*(\zeta(u))\leq \zeta^*(\eta(v)).$$
Going through the same argument but exchanging $\zeta$ and $\eta$ proves that $M$ is $\E$-principally symmetric.
\par
To prove the converse, we assume that $M$ is $\E$-principally generated. We showed in [Heymans and Stubbe, 2009, Prop.\ 8.2] that then necessarily $M$ is a locale and that there is a locale morphism\footnote{That locale morphism satisfies some further particular properties, which made us call it a {\em skew local homeomorphism} in that reference.} $f\:M\to X$ such that $m\tensor x=m\wedge f^*(x)$ for all $m\in M$ and $x\in X$. It follows easily from this characterisation that, for all $s\in\Gamma\can$, 
$$M\to M\:m\mapsto s\wedge m$$
is an $X$-module morphism. But under the hypothesis that $M$ is $\E$-principally symmetric, we can compose the left adjoint module morphism $\down\<s,s\>\can\to M\:x\mapsto s\tensor x$ with its right adjoint $M\to\down\<s,s\>\can\:m\mapsto \<s,m\>\can$ to obtain 
$$M\to M\:m\mapsto s\tensor\<s,m\>\can.$$
We claim that these module morphisms are equal: we shall show that they coincide on elements of $\Gamma\can$, which suffices because $\Gamma\can$ is a Hilbert basis. Indeed, for $r,t\in\Gamma\can$ we can compute that
$$\<r,s\wedge t\>\can=\<r,s\>\can\wedge\<r,t\>\can=\<r,s\>\can\wedge\<s,t\>\can=\<r,s\tensor\<s,t\>\can\>\can.$$
(The first equality holds because $\<r,-\>\can$ is a right adjoint, and the second equality holds because $\<r,s\>\can\wedge\<r,t\>\can=\<s,r\>\can\wedge\<r,t\>\can=\<s,r\tensor\<r,t\>\can\>\can\leq\<s,t\>\can$ and similarly $\<r,s\>\can\wedge\<s,t\>\can\leq\<r,t\>\can$.) Taking the supremum over all $r\in\Gamma\can$ proves that
$$s\wedge t=\bigvee_{r\in\Gamma\can}r\tensor\<r,s\wedge t\>\can=\bigvee_{r\in\Gamma\can}r\tensor\<r,s\tensor\<s,t\>\can\>\can=s\tensor\<s,t\>\can$$
as claimed. For any left adjoint $\zeta\:\down u\to M$ in $\Mod(X)$ we can apply the above to $s:=\zeta(u)\in\Gamma\can$, to find that 
$$\zeta(\zeta^*(m))=\zeta(u)\wedge m$$
for any $m\in M$. This allows us to verify in turn that for any $x\leq u$,
$$\zeta(x\wedge\zeta^*(m))=\zeta(\zeta^*(m))\tensor x=(\zeta(u)\wedge m)\tensor x=(\zeta(u)\tensor x)\wedge m=\zeta(x)\wedge m$$
as wanted.
\endofproof
The above direct argument relies on elementary order theory. There is a shorter alternative, using results in the literature: an \'etale $X$-module is the same thing as local homeomorphism into $X$ [Heymans and Stubbe, 2009, Theorem 7.12], which is the same thing as a Hilbert $X$-module [Rodrigues and Resende, 2008, Theorem 3.15], which is the same thing as an $\E$-principally generated and $\E$-principally symmetric $X$-module (by our upcoming Theorem \ref{14}).

\section{(Sometimes) all Hilbert structure is canonical}\label{D}

The previous section was concerned with the {\em canonical} Hilbert structure on a $Q$-module $M$: we showed that there is a {\em canonical} Hilbert basis for the {\em canonical} (pre-)inner product on $M$ if and only if $M$ is $\E$-principally generated and $\E$-principally symmetric, two natural notions based on the behaviour of certain adjunctions in $\Mod(Q)$. This section is devoted to the perhaps surprising fact that, for a certain class of quantales (containing many cases of interest), the {\em only possible} Hilbert structure is the canonical one. 
\begin{theorem}\label{14}
Let $Q$ be a modular quantal frame, $\E\subseteq Q$ the set of symmetric idempotents, and $M$ a $Q$-module. If $M$ bears a Hilbert structure, then necessarily $M$ is $\E$-principally generated and $\E$-symmetric, and the involved inner product is the canonical one, which moreover is strict (and, by Theorem \ref{11}, admits the canonical Hilbert basis).
\end{theorem}

The proof of the theorem shall be given as a series of lemmas. The first one straightforwardly extrapolates a result known to [Resende and Rodrigues, 2009] in the case of modules on a locale, and appears in [Resende, 2008]. We recalled the construction of the category $\Matr(Q)$ of matrices with entries in $Q$ in Example \ref{12}--\ref{12.4}, and remarked that whenever $Q$ is involutive then so is $\Matr(Q)$.
\begin{lemma}\label{15}
If $Q$ is an involutive quantale and $M$ is a $Q$-module with an inner product $\<-,-\>$ admitting a Hilbert basis $\Gamma$, then the following holds for all $m,n\in M$:
$$\bigvee_{s\in\Gamma}\<m,s\>\circ\<s,n\>=\<m,n\>.$$
In particular, $(\Gamma,\<-,-\>)$ is a so-called {\em projection matrix}: a symmetric idempotent in the involutive quantaloid $\Matr(Q)$.
\end{lemma}
\proof In $\<m,n\>$, use $n=\bigvee_{s\in\Gamma}s\tensor\<s,n\>$ and apply the linearity of $\<m,-\>$.
\endofproof

The following lemma refers to the notion of {\em total regularity}, which we here state in a bare-bones matrix-form, but which actually has deep connections with sheaf theory; it was introduced in the context of quantaloid-enriched categorical structures by Stubbe [2005b] and played a crucial role in [Heymans and Stubbe, 2009] too.
\begin{lemma}\label{16} Let $Q$ be an involutive quantale, and $M$ a $Q$-module with an inner product $\<-,-\>$ admitting a Hilbert basis $\Gamma$. The following statements are equivalent: 
\begin{enumerate}
\item\label{16.0} for all $s\in\Gamma$: $s=s\tensor\<s,s\>$,
\item\label{16.1} for all $s\in\Gamma$: $s\leq s\tensor\<s,s\>$,
\item\label{16.2} the projection matrix $(\Gamma,\<-,-\>)$ is {\em totally regular}, i.e.\ for all $s,t\in\Gamma$: 
$$\<s,t\>\circ\<t,t\>=\<s,t\>=\<s,s\>\circ\<s,t\>.$$
\end{enumerate}
If $Q$ moreover satisfies $q\leq q\circ q\o\circ q$ for every $q\in Q$, then these equivalent conditions always hold\footnote{More generally, under this condition it is true that any projection matrix with entries in $Q$, i.e.\ any symmetric idempotent in $\Matr(Q)$, is totally regular.}.
\end{lemma}
\proof 
Due to the Hilbert basis, $s\tensor\<s,s\>\leq\bigvee_{t\in\Gamma}t\tensor\<t,s\>=s$ for any $s\in\Gamma$ and thus $(\ref{16.1}\impl\ref{16.0})$. To see that ($\ref{16.1}\impl\ref{16.2}$), compute for $s,t\in\Gamma$ that $\<s,t\>=\<s,t\tensor\<t,t\>\>=\<s,t\>\circ\<t,t\>$. Conversely, ($\ref{16.2}\impl\ref{16.1}$) because, fixing a $t\in\Gamma$ we have for all $s\in\Gamma$ that $\<s,t\>=\<s,t\>\circ\<t,t\>=\<s,t\tensor\<t,t\>\>$; but therefore 
$$t=\bigvee_{s\in\Gamma}s\tensor\<s,t\>=\bigvee_{s\in\Gamma}s\tensor\<s,t\tensor\<t,t\>\>=t\tensor\<t,t\>.$$
Now if moreover every element $q\in Q$ satisfies $q\leq q\circ q\o\circ q$ then we can compute for $s,t\in\Gamma$ that
$$\<s,t\>\leq\<s,t\>\circ\<s,t\>\o\circ\<s,t\>=\<s,t\>\circ\<t,s\>\circ\<s,t\>
\left\{\begin{array}{c}\leq\<s,s\>\circ\<s,t\>\leq\<s,t\> \\ \mbox{ } \\ \leq\<s,t\>\circ\<t,t\>\leq\<s,t\> \end{array}\right.
$$
precisely as wanted in the second condition. (We used that $\<r,s\>\circ\<s,t\>\leq\<r,t\>$ for any $r,s,t\in\Gamma$, as follows trivially from the formula in Lemma \ref{15}.)
\endofproof
Particularly for a modular quantale $Q$ the above result is interesting: because $q\leq q\circ q\o\circ q$ holds as consequence of the modular law, it follows that for every $Q$-module with Hilbert structure its Hilbert basis is totally regular.

Next are two lemmas which contain the important (and less straightforward) matter.
\begin{lemma}\label{17}
If $Q$ is an involutive quantale, and $M$ a $Q$-module with an inner product $\<-,-\>$ admitting a Hilbert basis $\Gamma$ satisfying the equivalent conditions in Lemma \ref{16}, then for any $s\in\Gamma$ there is an adjunction
$$Q^{\<s,s\>}\xymatrix@=15mm{\ar@{}[r]|{\bot}\ar@/^2mm/@<1mm>[r]^{s\tensor-} & \ar@/^2mm/@<1mm>[l]^{\<s,-\>}}M$$
in $\Mod(Q)$. Writing $\E\subseteq Q$ for the set of symmetric idempotents, such an an $M$ is always $\E$-principally generated; and if $M$ is moreover $\E$-principally symmetric then $\<-,-\>$ coincides with the canonical (pre-)inner product $\<-,-\>\can$.
\end{lemma}
\proof 
For $s\in\Gamma$, compose the inclusion $Q^{\<s,s\>}\to Q$ with the module morphism $s\tensor-\:Q\to M$ to obtain a module morphism 
$$\zeta_s\,\:Q^{\<s,s\>}\to M:q\mapsto s\tensor q.$$ 
Because we assume $s=s\tensor\<s,s\>$ it follows that $\<s,s\>\circ\<s,m\>=\<s\tensor\<s,s\>,m\>=\<s,m\>$ for any $m\in M$, and therefore the obvious module morphism $\<s,-\>\:M\to Q$ co-restricts to a module morphism 
$$\zeta'_s\,\:M\to Q^{\<s,s\>}:m\mapsto\<s,m\>.$$
We shall show that $\zeta_s\dashv\zeta'_s$ in $\Mod(Q)$; in fact, it suffices to prove that this adjunction holds in the category of ordered sets and order-preserving maps. Thus, consider $q\in Q^{\<s,s\>}$ and $m\in M$: if $s\tensor q\leq m$ then $q=\<s,s\>\circ q=\<s,s\tensor q\>\leq\<s,m\>$; conversely, if $q\leq\<s,m\>$ then $s\tensor q\leq s\tensor\<s,m\>\leq\bigvee_{t\in\Gamma}t\tensor\<t,m\>=m$.

The module $M$ is $\E$-principally generated because for any $m\in M$ we have
\begin{eqnarray*}
m & = & \bigvee_{s\in\Gamma}s\tensor\<s,m\> \\
  & = & \bigvee_{s\in\Gamma}\zeta_s\zeta^*_s(m) \\
  & \leq & \bigvee\left\{\zeta(\zeta^*(m))\bigmid e\in\E,\ \zeta\:Q^e\to M\mbox{ left adjoint in }\Mod(Q)\right\} \\
  & \leq & m.
\end{eqnarray*}
It follows directly from Lemma \ref{15} that
$$\<m,n\>=\bigvee_{s\in\Gamma}(\zeta_s^*(m))\o\circ\zeta_s^*(n),$$
and from the above it is clear that this is smaller than
$$\<m,n\>\can=\bigvee\left\{(\zeta^*(m))\o\circ\zeta^*(n)\bigmid e\in\E, \zeta\:Q^e\to M\mbox{ left adjoint in }\Mod(Q)\right\}.$$
Now suppose that $M$ is $\E$-principally symmetric. Fixing a left adjoint $\zeta\:Q^e\to M$ in $\Mod(Q)$, with $e\in\E$, we can compute for any $s\in\Gamma$ that
$$\<\zeta(e),s\>=\<s,\zeta(e)\>\o=(\zeta^*_s(\zeta(e)))\o=\zeta^*(\zeta_s(\<s,s\>))=\zeta^*(s\tensor\<s,s\>)=\zeta^*(s);$$
the symmetry was crucially used in the third equality, and the assumption that the equivalent conditions in Lemma \ref{16} hold in the last one. But morphisms in $\Mod(Q)$ with domain $M$ are equal if they coincide on the Hilbert basis $\Gamma$, so for all $m\in M$ we have $\<\zeta(e),m\>=\zeta^*(m)$. Therefore we find that
\begin{eqnarray*}
(\zeta^*(m))\o\circ\zeta^*(n)
 & = & \<m,\zeta(e)\>\circ\zeta^*(n)\\
 & = & \<m,\zeta(e)\tensor\zeta^*(n)\>\\
 & = & \<m,\zeta(\zeta^*(n))\>\\
 & \leq &  \<m,n\>.
\end{eqnarray*}
To pass from the third to the fourth line we use the counit of the adjunction $\zeta\dashv\zeta^*$. All this means that $\<m,n\>\can\leq\<m,n\>$, and we are done.
\endofproof

It is only in the next statement that we require $Q$ to be a modular quantal frame.
\begin{lemma}\label{18}
If $Q$ is a modular quantal frame, and $M$ a $Q$-module with an inner product $\<-,-\>$ admitting a Hilbert basis $\Gamma$, then $M$ is $\E$-principally symmetric (for $\E\subseteq Q$ the set of symmetric idempotents). 
\end{lemma}
\proof 
We shall prove that the first of the equivalent conditions in Proposition \ref{10} holds {\em for the given inner product} on $M$; as we have remarked right after the proof of that Proposition, this suffices to infer the $\E$-principal symmetry of $M$. Because here we assume $M$ to have a Hilbert basis $\Gamma$, it in fact suffices to show that, for any $e\in\E$, any left adjoint $\zeta\:Q^e\to M$ and any $s\in\Gamma$: $\zeta^*(s)=\<\zeta(e),s\>$.

First remark that, with these notations, 
$$\zeta(e)\tensor\zeta^*(s)=\zeta(\zeta^*(s))\leq s$$
trivially holds. On the other hand, using all assumptions we can compute that
$$\begin{array}{llll}
e & = & e\wedge\zeta^*(\zeta(e)) & \mbox{(unit of $\zeta\dashv\zeta^*$)} \\[4pt]
 & = & \displaystyle{e\wedge\zeta^*\Big(\bigvee_{s\in\Gamma}s\tensor\<s,\zeta(e)\>\Big)} & \mbox{($\Gamma$ is a Hilbert basis)} \\[4pt]
 & = & \displaystyle{e\wedge\bigvee_{s\in\Gamma}\Big(\zeta^*\left(s\right)\circ\<s,\zeta(e)\>\Big)} & \mbox{($\zeta^*$ is a module morphism)} \\[4pt]
 & = & \displaystyle{\bigvee_{s\in\Gamma}\Big(e\wedge\zeta^*\left(s\right)\circ\<s,\zeta(e)\>\Big)} & \mbox{($Q$ is a frame)}\\[4pt]
 & \leq & \displaystyle{\bigvee_{s\in\Gamma}\Big(e\circ\<s,\zeta(e)\>\o\wedge\zeta^*\left(s\right)\Big)\circ\<s,\zeta(e)\>} & \mbox{(by the modular law)}\\[4pt]
 & = & \displaystyle{\bigvee_{s\in\Gamma}\Big(\<\zeta(e),s\>\wedge\zeta^*\left(s\right)\Big)\circ\<s,\zeta(e)\>} & \mbox{(symmetry)} \\[4pt]			 
 & \leq	& \displaystyle{\bigvee_{s\in\Gamma}\<\zeta(e),s\>\circ\<s,\zeta(e)\>} & \mbox{(trivially)} \\[4pt]
 & = & \displaystyle{\<\zeta(e),\zeta(e)\>}. & \mbox{(by Lemma \ref{15})}			 
\end{array}$$
Hence, combining both the previous inequalitites,
$$\zeta^*(s)=e\circ\zeta^*(s)\leq\<\zeta(e),\zeta(e)\>\circ\zeta^*(s)=\<\zeta(e),\zeta(e)\tensor\zeta^*(s)\>\leq\<\zeta(e),s\>$$
and we have the ``$\leq$'' of the required equality. To see that also ``$\geq$'' holds, we first apply the modularity of $Q$ again to compute
$$\begin{array}{llll}
e & = & \displaystyle{\bigvee_{s\in\Gamma}\Big(e\wedge\zeta^*\left(s\right)\circ\<s,\zeta(e)\>\Big)} & \mbox{(as above)} \\[4pt]
  & \leq & \displaystyle{\bigvee_{s\in\Gamma}\zeta^*\left(s\right)\circ\Big(\left(\zeta^*\left(s\right)\right)\o\circ e\wedge\<s,\zeta(e)\>\Big)} & \mbox{(by the modular law)}\\[4pt]
  & \leq & \displaystyle{\bigvee_{s\in\Gamma}\zeta^*\left(s\right)\circ\left(\zeta^*\left(s\right)\right)\o}. & \mbox{(trivial)}  
\end{array}$$
Now we combine this with the first inequality that we proved to obtain
$$\begin{array}{llll}
\<\zeta(e),s\>
 & = & e\circ\<\zeta(e),s\> & \mbox{(trivial)} \\[4pt]
 & \leq & \displaystyle{\bigvee_{t\in\Gamma}\zeta^*\left(t\right)\circ\left(\zeta^*\left(t\right)\right)\o\circ\<\zeta(e),s\>} & \mbox{(by the above)} \\[4pt]
 & = & \displaystyle{\bigvee_{t\in\Gamma}\zeta^*\left(t\right)\circ\<\zeta(e)\tensor\zeta^*\left(t\right),s\>} & \mbox{(``conjugate-linearity'' of inner product)} \\[4pt]
 & \leq & \displaystyle{\bigvee_{t\in\Gamma}\zeta^*\left(t\right)\circ\<t,s\>} & \mbox{(because $\zeta(e)\tensor\zeta^*(t)\leq t$)} \\[4pt]	
 & = & \displaystyle{\zeta^*\Big(\bigvee_{t\in\Gamma}t\tensor\<t,s\>\Big)} & \mbox{($\zeta^*$ is a module morphism)} \\[4pt]
 & = & \displaystyle{\zeta^*(s)}. &	\mbox{($\Gamma$ is a Hilbert basis)} 
\end{array}$$
and we are done.
\endofproof
Relying on categorical machinery, the previous Lemma can alternatively be proved as follows: Bearing in mind Example \ref{12}--\ref{12.5}, the requirement that $\zeta^*(s)=\<\zeta(u),s\>$ at the end of the first paragraph of the proof above is equivalent to asking for $\zeta^*=\zeta\o$ in the category $\Hilb(Q)$. But $Q$ being a modular quantal frame is equivalent to the matrix quantaloid $\Matr(Q)$ being modular, which in turn implies that the quantaloid $\Proj(Q)$ of projection matrices, obtained by splitting the symmetric idempotents in $\Matr(Q)$, is modular too. In any modular quantaloid, the right adjoint of a morphism, should it exist, is necessarily its involute: this thus holds in $\Proj(Q)$, and also in its equivalent $\Hilb(Q)$. Therefore in particular $\zeta^*=\zeta\o$ for any left adjoint $\zeta\:Q^e\to M$, as wanted.

Lastly we have a simple lemma asserting the strictness of inner products in certain cases.
\begin{lemma}\label{19}
Let $Q$ be an involutive quantale in which $q\leq q\circ q\o\circ q$ holds for any $q\in Q$. If $M$ is a $Q$-module with an inner product $\<-,-\>$ admitting a Hilbert basis $\Gamma$, then this inner product is strict.
\end{lemma}
\proof Let $q\in Q$: if $q\o\circ q=0$ then $q\leq q\circ q\o\circ q=q\circ0=0$ hence $q=0$. Now suppose that $\<m,m\>=0$ for an $m\in M$. The formula in Lemma \ref{15} implies that 
$$\<s,m\>\o\circ\<s,m\>=\<m,s\>\circ\<s,m\>\leq\bigvee_{t\in\Gamma}\<m,t\>\circ\<t,m\>=0$$ 
for all $s\in\Gamma$, whence $\<s,m\>=0$ for all $s\in\Gamma$. But then $m=\bigvee_{s\in\Gamma}s\tensor\<s,m\>=0$ as required.
\endofproof

Having all these lemmas, we assemble the proof of the statement in the beginning of this section.
\par\medskip\noindent {\em Proof of Theorem \ref{14}\ }: 
Because $Q$ is by hypothesis a modular quantal frame we have by Lemma \ref{18} that $M$ is $\E$-principally symmetric. It follows from $Q$'s modularity and Lemma \ref{16} that Lemma \ref{17} applies, showing that $M$ is $\E$-principally generated. Together with the fact that $M$ is $\E$-principally symmetric this moreover entails the equality of the given inproduct with the canonical one. Finally, the strictness of the (canonical) inner product is a consequence of Lemma \ref{19}.
\endofproof

\begin{example}\label{20} We end with examples that refer to the category $\Hilb(Q)$ of Hilbert modules, and particularly to applications in sheaf theory.
\begin{enumerate}
\item\label{20.1} As in Example \ref{12}--\ref{12.4} we write $\Hilb(Q)$ for the quantaloid of $Q$-modules with Hilbert structure. For a modular quantal frame $Q$, Theorem \ref{14} allows us to consider $\Hilb(Q)$ as a full subquantaloid of $\Mod(Q)$: there is only one relevant Hilbert structure on a $Q$-module. Moreover, Lemma \ref{18} implies that, whenever $\phi\:M\to N$ is a left adjoint in $\Hilb(Q)$, then $\phi\dashv\phi\o$ (compare with Example \ref{12}--\ref{12.5}). Because in this case every symmetric idempotent in $\Matr(Q)$ is totally regular, we therefore get the equivalences of quantaloids
$$\Hilb(Q)\simeq\Proj(Q)\simeq\Dist\o(Q_\E)$$
where $Q_\E$ denotes the quantaloid obtained as universal splitting of the symmetric idempotents of $Q$, and $\Dist\o(Q_\E)$ is the full subquantaloid of $\Dist(Q_\E)$ (= the quantaloid of $Q_\E$-enriched categories and distributors [Stubbe, 2005a]) determined by the {\em symmetric $Q_\E$-categories}.  
\item\label{20.2} Sheaves on sites: For a small site $(\C,J)$ and $Q$ the associated modular quantal frame as in Example \ref{1.1}--\ref{1.1.6}, the category $\Sh(\C,J)$ is equivalent to the category of $Q$-modules with canonical Hilbert structure and the left adjoint module morphisms between them:
$$\Sh(\C,J)\simeq\Map(\Hilb(Q)).$$
With a bit more work this can be rephrased as equivalent quantaloids: 
$$\Rel(\Sh(\C,J))\simeq\Dist\o(\Q)\simeq\Hilb(Q).$$
\par\noindent\textit{Sketch of the proof:} Let $\Q$ be as in Example \ref{1.1}--\ref{1.1.6}. Walters [1982] proved that $\Sh(\C,J)$ is equivalent to $\Map(\Dist\o(\Q))$, the full subcategory of $\Map(\Dist(\Q))$ determined by the \textit{symmetric} $\Q$-categories. In the previous example we indicated that $\Hilb(Q)\simeq\mathsf{Proj}(Q)\simeq\Dist\o(Q_\E)$, hence it suffices to prove that $\Dist(\Q)\simeq\Dist(Q_\E)$. But this follows from the fact that $\Q$, regarded as a subquantaloid of $Q_\E$, is {\em dense} in $Q_\E$: for any $X\in Q_\E$, $\id_X=\bigvee_{i\in I}f_i\circ f_i^*$ with each $f_i\:X_i\to X$ a left adjoint with $X_i\in\Q_0$. This property is due to the modularity of $Q$ and the coreflexive idempotents ($e\leq\id_{\dom(e)}$) of $\Q$ (corresponding to closed sieves on $\C$) being suprema of the form above.
\endofproof	
\item\label{20.3} Sheaves on an \'etale groupoid: We understand that Resende [2008] defines a ``sheaf'' on an involutive quantale $Q$ to be a $Q$-module $M$ with Hilbert structure {\em satisfying moreover} $\bigvee\Gamma=\top_M$, and proves -- via the correspondence between \'etale groupoids and inverse quantal frames from [Resende, 2007] -- that, for an \'etale groupoid $\G$, the topos of $\G$-sheaves is equivalent to the category with as objects those ``sheaves'' on an inverse quantal frame $\O(\G)$ and as morphisms the left adjoint $\O(\G)$-module morphisms that have their involute as right adjoint (which he describes as ``direct image homomorphisms''). This may appear to be in contradiction with the examples above: sheaves (on a site) can be described as modules (on a modular quantal frame) with Hilbert structure {\em without} any further conditions. However it turns out that, when $Q$ is an inverse quantal frame (as in the main example of [Resende, 2008]), then the extra condition is anyway a {\em consequence} of the features of $Q$ (see the proof below). We further repeat from Example \ref{20}--\ref{20.1} that, because $\O(\G)$ is a modular quantal frame, any left adjoint in $\Hilb(\O(\G))$ has its involute as right adjoint (but this need not be so for involutive quantales in general, where we think this is an important extra condition). Conclusively, by Theorem \ref{14} the topos of sheaves on an \'etale groupoid $\G$ is equivalent to $\Map(\Hilb(\O(\G)))$.
\par\noindent\textit{Proof:} If $Q$ is an inverse quantal frame and $M$ a $Q$-module with inner product $\<-,-\>$ admitting a Hilbert basis $\Gamma$, we may assume without loss of generality that $\Gamma$ is {\em maximal} in the following way: 
$$\Gamma=\left\{s\in M\bigmid\forall m\in M:s\tensor\<s,m\>\leq m\right\}.$$
If $p\in Q$ is a partial unit and $s\in\Gamma$ then $s\tensor p\in\Gamma$: because for all $m\in M$, 
$$(s\tensor p)\tensor\<s\tensor p,m\>=(s\tensor pp\o)\tensor\<s,m\>\leq s\tensor\<s,m\>\leq m;$$
hence certainly $s\tensor p\leq\bigvee\Gamma$. Since $\top_Q$ is by assumption the join of all partial units, this implies $s\tensor\top_Q\leq\bigvee\Gamma$, whence 
$$\top_M=\bigvee_{s\in\Gamma}s\tensor\<s,\top_M\>\leq\bigvee_{s\in\Gamma}s\tensor\top_Q\leq\bigvee\Gamma$$
which proves the claim.
\endofproof
\end{enumerate}
\end{example}

\section{Concluding remarks}\label{E}

In this paper we proved the following results in the theory of quantale modules: (1) every module on an involutive quantale $Q$ bears a canonical (pre-)inner product; (2) that canonical (pre-)inner product admits the canonical Hilbert basis if and only if the module is principally generated and principally symmetric; and (3) if $Q$ is a modular quantal frame then the only possible Hilbert structure (= inner product plus Hilbert basis) on a $Q$-module is the canonical one. In the examples we explained the use of these results in sheaf theory: we argued in particular that the category of sheaves on a site $(\C,J)$ is equivalent to a category of quantale modules with (canonical) Hilbert structure.

These results are a natural continuation of our previous work. Whereas Stubbe [2005b] described ordered sheaves on a quantaloid $\Q$ as particular $\Q$-enriched categorical strutures, Heymans and Stubbe [2009] reformulated this -- via the correspondence between cocomplete $\Q$-categories and $\Q$-modules, and the particular role of $\Q$-modules in the theory of ordered sheaves on $\Q$ [Stubbe, 2006, 2007] -- in a module-theoretic language: ordered sheaves on $\Q$ are the same thing as principally generated $\Q$-modules. The material in this paper suggests that the ``symmetrically ordered'' sheaves (i.e.\ sheaves {\it tout court}) on an involutive quantale $Q$ are those principally generated $Q$-modules which are moreover principally symmetric. The latter in turn coincide with modules bearing a canonical Hilbert structure (which, for modules on a modular quantal frame, is the only possible Hilbert structure). 

Our current research is concerned with a further elaboration of that novel notion, ``principal symmetry'': we extend it from quantale modules to quantaloid modules, and even to quantaloid-enriched categories. A future paper shall in particular contain all remaining details from Examples \ref{1.1}--\ref{1.1.6} and \ref{20}--\ref{20.2}.

\end{document}